\documentclass[11pt]{article}

\usepackage{amsmath,amssymb}

\usepackage{color}

\renewcommand{\subsection}{\subsubsection}

\newtheorem{theorem}{Theorem}[section]

\newtheorem{remark}{Remark}[section]

\def\div{{\rm div}\,}
\def\diag{{\rm diag}\,}
\def\a{\alpha}

\begin{document}

\title{\bf Symmetrizations of RMHD equations and \\ stability of relativistic current-vortex sheets}

\author{{\bf Heinrich Freist\"uhler}\\
Universit\"at Konstanz, Fachbereich Mathematik und Statistik,\\ Fach 199, 78457 Konstanz, Germany\\
E-mail: Heinrich.Freistuehler@uni-konstanz.de
\and
{\bf Yuri Trakhinin}\\
Sobolev Institute of Mathematics, Koptyug av. 4, 630090 Novosibirsk, Russia\\
E-mail: trakhin@math.nsc.ru
}

\date{}

\maketitle

\begin{abstract}
We consider the equations of relativistic magnetohydrodynamics (RMHD) in the case of special relativity.
{Starting by computations in the fluid's rest frame and then applying Lorentz
transformations, we derive a covariant symmetric formulation of RMHD in terms of the primitive
(physical) variables. This symmetric system is important for the study of various
initial boundary value problems.} 
We also find a so-called
secondary symmetrization whose direct consequence is the extension of the sufficient
stability condition obtained earlier for non-relativistic planar current-vortex sheets
to the relativistic case. As in non-relativistic settings, this implies the local-in-time
existence of corresponding smooth nonplanar current-vortex sheets.
\end{abstract}

\section{Introduction}
\label{sec:1}
The mathematical model of magnetohydrodynamics has been very widely used in the literature
of the last decades to study self-consistent interactions of flows of charged matter with
magnetic fields in terrestrial and non-terrestrial contexts \cite{Schindler}.
In important astrophysical situations, velocities are known to become so large
that one cannot ignore the effects of relativity.
In such cases, e.\ g., of fast stellar winds or jets, or explosive plasma outflows
from collapsing stars,  one uses the special relativistic version RMHD of magnetohydrodyanmics.
We exemplarily refer to \cite{Kom1999} and \cite{Mat2011} both
regarding detailed references for such contexts and, in particular, modern numerical schemes
for RMHD.

The present, theoretically oriented paper contributes to the understanding of the RMHD equations
by deriving a family of new symmetric versions thereof. While
we do think that our new formulations of RMHD have a good chance to prove useful also
for computational purposes, their application in the current paper is a different one: These
symmetric versions allow for rigorous analytical investigation of discontinuous flow patterns, as are
common in magnetospheres and other flows of charged particles, regarding their dynamical stability!
We show this at the challenging example of relativistic current-vortex sheets.
\par\bigskip

The equations of relativistic magnetohydrodynamics (RMHD) for an ideal fluid read (see, e.g., \cite{An,Lich}):
\begin{equation}
\nabla_{\a}(\rho u^{\a})=0, \qquad
\nabla_{\a}T^{\a\beta}=0,
\qquad
\nabla_{\a}(u^{\a}b^{\beta} -u^{\beta}b^{\a})=0,
\label{1}
\end{equation}
where $\nabla_{\a}$ is the covariant derivative with respect to
the Lorentzian metric $g=\diag (-1,1,1,1)$ of the {space-time with 
components} $g_{\a\beta}$ (in this paper we restrict ourselves to the case of {\it special relativity}), $\rho$ is the proper rest-mass density of the plasma, $u^{\a}$ are components of the 4-velocity,
\[
T^{\a\beta}=(\rho h +B^2)u^{\a}u^{\beta} +qg^{\a\beta} -b^{\a}b^{\beta},
\]
$h= 1+e +(p/\rho )$ is the relativistic specific enthalpy, $p$ is the pressure, $e=e(\rho ,S)$ is the specific internal energy, $S$ is the specific entropy, $B^2=b^{\a}b_{\a}$, $b^{\a}$ are components of the magnetic field 4-vector with respect to the plasma velocity, and $q= p+\frac{1}{2}B^2$ is the total pressure. The 4-vectors satisfy the conditions
$u^{\a}u_{\a} =-1$ and $u^{\a}b_{\a} =0$, and the speed of the light is equal to unity.

Let $(x^0,x)$ be inertial coordinates, with $t=x^0$ a time coordinate and $x=(x^1,x^2,x^3)$ space coordinates. Then
\[
u^0=-u_0=\Gamma ,\quad u^i=u_i =\Gamma v_i,\quad u:=(u_1,u_2,u_3)=\Gamma v,
\]
\[
\Gamma^2=1+|u|^2,\quad
b^0=-b_0=(u,H),\quad
b^i=b_i=\frac{H_i}{\Gamma}+(u,H)v_i,
\]
\[
b:=(b_1,b_2,b_3),\quad
B^2=|b|^2-b_0^2=\frac{|H|^2}{\Gamma^2} +(v,H)^2>0,
\]
where $\Gamma=(1-|v|^2)^{-1/2}$ is the Lorentz factor, $v=(v_1,v_2,v_3)$ is the plasma 3-velocity, $H=(H_1,H_2,H_3)$ is the magnetic field 3-vector, the Latin indices run from 1 to 3, and by $(\ ,\ )$ we denote the scalar product.

The RMHD equations \eqref{1} can be written as the system of conservation laws
\begin{align}
& \partial_t(\rho\Gamma ) +\div (\rho u )=0,\label{2}\\
& \partial_t\left(\rho h \Gamma u +|H|^2v-(v,H)H\right) \nonumber \\
&  + {\rm div} \left((\rho h +B^2) u\otimes u  -b\otimes b\right) + {\nabla}q=0,\label{3}\\
& \partial_t(\rho h\Gamma^2 +|H|^2-q ) +{\rm div} \left(\rho h\Gamma u +|H|^2v-(v,H)H \right)=0,\label{4}\\
& \partial_tH -{\nabla}\times (v {\times} H)=0,\label{5}
\end{align}
where $\partial_t=\partial /\partial t$, $\nabla =(\partial_1,\partial_2,\partial_3)$, $\partial_i=\partial /\partial x^i$,
etc., and the constraint 
\begin{equation}
\div H=0.\label{6}
\end{equation}
The state representation $U=(p,u,H,S)$ being known as a natural choice of {\it primitive} (physical) variables \cite{Font},
we will take $p,u,H,S$ as primary unknowns.

The first two main results of this paper can be formulated as follows.

\begin{theorem} For smooth solutions satisfying the constraint \eqref{6},
the RMHD equations \eqref{2}--\eqref{5}  are equivalently rewritten as the symmetric system
\begin{equation}
A_0(U)\partial_tU+\sum_{j=1}^3 A_j(U)\partial_jU=0
\label{first}
\end{equation}
for the vector $U=(p,u,H,S)$, provided that the hyperbolicity condition $A_0>0$ holds, i.\ e,.
\begin{equation}
\rho >0,\quad p_{\rho} >0,\quad 0<c_s^2<1,\label{hc}
\end{equation}
where $c_s^2= a^2/h$ is the relativistic speed of sound, $a^2=p_{\rho}(\rho ,S)$,
\begin{equation}
A_0=
\begin{pmatrix}
{\displaystyle\frac{\Gamma}{\rho a^2}} & v^{\sf T} &0 &0 \\[3pt]
v & \mathcal{A} &0 &0 \\
0& 0 & \mathcal{M} &0 \\
0 &0&0 & 1
\end{pmatrix},\qquad A_j=
\begin{pmatrix}
{\displaystyle\frac{u_j}{\rho a^2}} & e_j^{\sf T} &0 &0\\[3pt]
e_j & \mathcal{A}_j &{\mathcal{N}_j}^{\sf T} &0\\
0& \mathcal{N}_j & v_j\mathcal{M} &0 \\
0&0&0&v_j
\end{pmatrix},\label{27}
\end{equation}
$e_j= (\delta_{1j},\delta_{2j},\delta_{3j})$ are the unit column vectors,
\[
\mathcal{A}=\left( \rho h \Gamma +\frac{|H|^2}{\Gamma} \right) I - \left( \rho h \Gamma +\frac{|H|^2+B^2}{\Gamma}\right) v\otimes v   -\frac{1}{\Gamma}\, H\otimes H +
\frac{(v,H)}{\Gamma}
\bigl( v\otimes H + H\otimes v\bigr) ,
\]
\[
\mathcal{M} =\frac{1}{\Gamma}\, (I +u\otimes u),\qquad \mathcal{N}_j= \frac{1}{\Gamma}\,b\otimes e_j -\frac{v_j}{\Gamma}\, b\otimes v - \frac{H_j}{\Gamma^2}\,I,
\]
\begin{multline*}
\mathcal{A}_j=
v_j\left\{ \left( \rho h \Gamma +\frac{|H|^2}{\Gamma} \right) I - \left( \rho h \Gamma +\frac{|H|^2-B^2}{\Gamma}\right) v\otimes v  -\frac{1}{\Gamma}\, H\otimes H\right\}\\
+\frac{H_j}{\Gamma}\left\{ \frac{1}{\Gamma^2}
\bigl( v\otimes H + H\otimes v\bigr) -2(v,H)(I-v\otimes v) \right\}
\\
+\frac{(v,H)}{\Gamma}\left( H\otimes e_j + e_j \otimes H\right) -\frac{B^2}{\Gamma}
\left( v\otimes e_j + e_j \otimes v\right),
\end{multline*}
and $I$ is the unit matrix.
\label{t1}
\end{theorem}

\begin{theorem}
More generally, there is a non-trivial one-parameter family of matrix-field quadruples
$(A_0^\lambda,A_1^\lambda,A_2^\lambda,A_3^\lambda)$ with
$$
A_\alpha^\lambda(U)\equiv A_\alpha(U) +\lambda \hat A_\alpha(U),\quad\alpha=0, \ldots,3,\qquad \lambda\in
\mathbb{R},
$$
such that for each $\lambda$,
the RMHD equations \eqref{2}--\eqref{5}, again assuming (6),
are equivalently rewritten as a symmetric system
\begin{equation}
A_0^\lambda(U)\partial_tU+\sum_{j=1}^3 A_j^\lambda(U)\partial_jU=0,
\label{second}
\end{equation}
provided that the hyperbolicity condition $A_0^\lambda >0$ holds, i.\ e., conditions \eqref{hc} are satisfied together with the inequality
\begin{equation}
\label{hyprest}
\lambda^2<\frac{c_s^2}{\rho  a^2 + B^2}.
\end{equation}
The symmetric matrices $A_{\alpha}^\lambda$ are specified in Section 3.
\label{t2}
\end{theorem}

A time-dependent hypersurface $\Sigma(t)\subset\mathbb{R}^3$
is called a current-vortex sheet if the velocity and magnetic-field components $v_N^\pm, H_N^\pm$
normal to $\Sigma$ are
always $0$ and $q$ does not jump across $\Sigma$ (while jumps in the tangential components
of $v$ and $H$ indicate vorticity and electric current along $\Sigma$).
\par\bigskip
The following is the third main result of this article.

\begin{theorem}
If the unperturbed uniform (piecewise constant) flow for a planar relativistic
current-vortex sheet satisfies a certain natural open condition (\eqref{sc}, with the notations in
\eqref{not}) and the unperturbed tangential magnetic fields on either side
of the discontinuity are nonzero and nonparallel to each other, then this current-vortex
sheet is linearly stable.
\label{t3}
\end{theorem}

Theorems \ref{t1}, \ref{t2}, \ref{t3} will be shown in Sections 2, 3, and 4, respectively. In the remaining part
of this introduction we discuss the contents of the theorems and their relations with each other
as well as with previously existing knowledge.
\par\bigskip
The fact that RMHD can at all be written symmetrically also in
the relativistic case has been known for long time.
Using \eqref{6} and the additional conservation law (entropy conservation)
\begin{equation}
\partial_t(\rho\Gamma S ) +\div (\rho S u )=0\label{7},
\end{equation}
which holds on smooth solutions of system \eqref{2}--\eqref{5}, and following Godunov's
symmetrization procedure \cite{Go1,Go} (see also, e.g., \cite{BR,Bl,BT} for its description)
one can symmetrize the conservation laws \eqref{2}--\eqref{5} in terms of a vector of canonical
variables $Q=Q(U)$. This was pointed out by Ruggeri and Strumia \cite{RuSt} and also by
Anile and Pennisi \cite{AP,An}. A concrete form of symmetric matrices was, however, not given
in \cite{RuSt,AP}. Also, dealing with initial-boundary value problems, notably with
free boundaries, it would be very inconvenient to work in terms of the vector $Q$.
In principle (and once one had a concrete from of the matrices with respect to the $Q$-representation)
one could transform from $Q$ to another set of primitive variables, in particular $U$,
keeping the symmetry property (see \cite{BT}).
However, finding a concrete form of symmetric matrices associated to the vector $Q$
and then doing that transformation from $Q$ to $U$ seems
connected with unimaginably lengthy calculations.

By contrast, our symmetric systems \eqref{first}, \eqref{second} in Theorems \ref{t1}, \ref{t2}
are nothing else than {\it algebraic symmetrizations} of the RMHD equations. That is,
following Friedrichs and Lax \cite{FL}, one first rewrites a given system of conservation laws in
the quasilinear form
\begin{equation}
B_0(U)\partial_tU+\sum_{j=1}^3 B_j(U)\partial_jU=0
\label{quasi}
\end{equation}
and then tries to somehow find a nonsingular matrix $D(U)$ such that
multiplication of \eqref{quasi} from the left by this matrix,
\[
D(U)B_0(U)\partial_tU+\sum_{j=1}^3 D(U)B_j(U)\partial_jU=A_0(U)\partial_tU+\sum_{j=1}^3
A_j(U)\partial_jU=0,
\]
makes 
\begin{equation}
A_0=A_0^{\sf T} > 0 \quad \mbox{and}\quad A_j=A_j^{\sf T}
\label{symm.cond}
\end{equation}
(where the hyperbolicity condition $A_0>0$ may sometimes require restriction
of the unknown $U$ to a subset of its original range).
This simple procedure of algebraic symmetrization was proposed many years ago by Friedrichs
and the matrix $D$ is sometimes called {\it Friedrichs' symmetrizer}.

However, we should modify the above procedure for the case when the system of conservation laws is supplemented by a set of divergence constraints on the initial data:
\begin{equation}
\div\Psi_i(U)=0,\quad i=1,\ldots,K.
\label{divcon}
\end{equation}
Namely, we do not just multiply \eqref{quasi} from the left by $D(U)$ but also take into account constraints \eqref{divcon}:
\[
D(U)B_0(U)\partial_tU+\sum_{j=1}^3 D(U)B_j(U)\partial_jU+\sum_{i=1}^KR_i(U)\,\div\Psi_i(U)
\]
\begin{equation}
=A_0(U)\partial_tU+\sum_{j=1}^3 A_j(U)\partial_jU=0,
\label{symmet}
\end{equation}
where $R_i(U)$ are some vector-columns. It is natural to call the set $\{ D, R_1,\ldots, R_K\}$
the generalized Friedrichs' symmetrizer if system \eqref{symmet} satisfies conditions
\eqref{symm.cond}. Thus, as in Godunov's symmetrization procedure \cite{Go,BR,Bl,BT},
constraints play an important role in the process of algebraic symmetrization.

For the RMHD equations we have only one divergence constraint, \eqref{6}. Taking into account \eqref{6}, in Sect. \ref{sec:2} we first perform the described algebraic symmetrization of the RMHD system for the fluid rest frame ($v=0$) and then apply the Lorentz transformation to the resulting symmetric system.
In fact, in Sect.\ \ref{sec:2} we do not write down the symmetrizer $\{ D, R_1\}$ concretely
but just obtain a symmetric form \eqref{symmet} of the RMHD system for the fluid rest frame by directly
rewriting it,
in view of \eqref{6}, in a nonconservative form for a suitable choice of the unknown $U$.
As for the non-relativistic MHD \cite{BT}, the good choice is  $U=(p,u,H,S)$.


We briefly compare the situation with classical fluid dynamics.
In \cite{Go1}, Godunov first symmetrized the compressible Euler equations in terms of the vector of canonical variables $Q=(q_1,\ldots,q_5)$, with
\[
q_1= -\frac{1}{T}\left( e +pV-ST-\frac{|v|^2}{2}\right),\quad q_{1+k}=-\frac{v_k}{T},\quad k=1,2,3,\quad q_5=\frac{1}{T},
\]
where $T$ is the temperature, $e=e(V,S)$ is the internal energy, and $V=1/\rho$. The corresponding symmetric system (for its concrete form we refer to \cite{Go1} or \cite{Bl}) is hyperbolic if in addition to the natural physical restrictions (cf. \eqref{hc})
\begin{equation}\label{hc1}
\rho >0,\quad a^2=p_{\rho} >0
\end{equation}
the convexity conditions for the function  $e=e(V,S)$ hold. At the same time, it is well-known that the compressible Euler equations rewritten in the nonconservative form
\begin{equation}
\label{gd}
\frac{1}{\rho a^2}\,\frac{{\rm d} p}{{\rm d}t} +{\rm div}\,{v} =0,\qquad
\rho\, {\displaystyle\frac{{\rm d}v}{{\rm d}t}+{\nabla}p  =0 },\qquad
{\displaystyle\frac{{\rm d} S}{{\rm d} t} =0}
\end{equation}
forms a symmetric system for the unknown $U=(p,v,S)$ which is hyperbolic under the physical restrictions \eqref{hc1}. The algebraic symmetrization \eqref{gd} and the Godunov's symmetrization in \cite{Go1} are equivalent if the more restrictive hyperbolicity condition for the symmetric form in  \cite{Go1} is satisfied,
i.\ e., if inequalities \eqref{hc1} hold together with the convexity conditions for $e=e(V,S)$.
Now,
roughly speaking, the relation between the algebraic symmetrization \eqref{first} of the RMHD equations and
their Godunov's symmetrization in \cite{RuSt,AP,An} is the same as the relation between the nonconservative
form \eqref{gd} of the Euler equations and their Godunov's symmetrization in \cite{Go1}.
In particular, \eqref{first} is equivalent to the symmetric system in \cite{RuSt,AP,An}, provided
that the physical restrictions \eqref{hc} are supplemented with the convexity conditions from
\cite{RuSt,BR}. We again underline that our main goal was to find a symmetric form of the RMHD equations
convenient for its usage for initial boundary value problems in regions of smoothness, and that in this
regard the knowledge of a {\it concrete}\ form of symmetric matrices in \eqref{first} is very important.

It is interesting to note that one system of conservation laws can have more than
one algebraic symmetrization. In Section 3, as in the non-relativistic case \cite{T05,T09},
we find --- besides the first symmetric version in Theorem \ref{t1} ---, as Theorem \ref{t2},
a so-called {\it secondary} symmetrization.
We call it secondary because the RMHD system being written in form \eqref{first} is already symmetric,
but we again find for \eqref{first} a nontrivial generalized Friedrichs' symmetrizer.
The hyperbolicity restriction \eqref{hyprest} is compatible with what we do with Theorem \ref{t2}
in the proof of Theorem \ref{t3}.

Finally we do treat an inital-boundary value problem in this paper.
Using our first (Theorem \ref{t1}) and second(ary) algebraic symmetrizations (Theorem \ref{t2})
of the RMHD equations, in Section 4 we derive a sufficient stability condition for planar relativistic
current-vortex sheets; this is Theorem \ref{t3}. As for the non-relativistic case \cite{T05,T09},
this stability condition being satisfied at each point of an initial nonplanar
current-vortex sheet guarantees the local-in-time existence of current-vortex sheet solutions
of the original nonlinear system. Here, the secondary symmetrization enables one to make the boundary
conditions of the linearized problem dissipative. It is actually the hyperbolicity condition
\eqref{hyprest} for the secondary symmetrization that gives us the sufficient stability condition
for planar relativistic current-vortex sheets.
\par\bigskip
Last not least, we point out that the first symmetrization \eqref{first} has already been crucially used
in \cite{T10} to prove the stability of the relativistic plasma-vacuum interface.

\section{Symmetrization of RMHD in terms of primitive variables}
\label{sec:2}

Note that equations \eqref{2} and \eqref{7} imply
\begin{equation}
\frac{{\rm d} S}{{\rm d} t} =0,
\label{8}
\end{equation}
with ${\rm d} /{\rm d} t =\partial_t+({v} ,{\nabla} )$. In the absence of the magnetic field, $H=0$, following \cite{T09.cpam}, from \eqref{2}--\eqref{4},  \eqref{8} we obtain the nonconservative form of the relativistic Euler equations
\begin{gather}
\frac{\Gamma}{\rho a^2}\,\frac{{\rm d} p}{{\rm d} t} +(v,\partial_tu) +{\rm div}\, u =0,\label{9}\\
(\rho h\Gamma )\left ( \frac{{\rm d} u}{{\rm d} t} -v \left( v,\frac{{\rm d} u}{{\rm d} t} \right)\right)
+(\partial_t p)v +\nabla p =0,\qquad
\frac{{\rm d} S}{{\rm d} t} =0.
\label{10}
\end{gather}
Equations \eqref{9}, \eqref{10} form the symmetric system
\begin{equation}
B_0(W)\partial_tW+\sum_{j=1}^3 B_j(W)\partial_jW=0
\label{11}
\end{equation}
for $W=(p,u,S)$, where the matrices $B_{\alpha}$ were written in \cite{T09.cpam}:
\[
B_0=\begin{pmatrix}
{\displaystyle\frac{\Gamma}{\rho a^2}} & v^{\sf T} &0  \\[3pt]
v & \rho h\Gamma\mathcal{B} &0  \\
0&0 & 1
\end{pmatrix},\qquad B_j=
\begin{pmatrix}
{\displaystyle\frac{u_j}{\rho a^2}} & e_j^{\sf T} &0 \\[3pt]
e_j & \rho hu_j\mathcal{B} &0 \\
0&0&v_j
\end{pmatrix},\qquad
\mathcal{B}=I-v\otimes v.
\]

The natural idea of symmetrizing the RMHD equations in the general case of a nonzero magnetic field is the repetition of simple arguments above that could give the desired symmetric system \eqref{first} for $U=(p,u,H,S)$, and it is reasonable to expect that the matrices $A_{\alpha}$ with cancelled fifth, sixth and seventh rows and columns coincide with $B_{\alpha}$ in \eqref{11} for $H=0$. Unfortunately, this way is connected with extremely long technical calculations. To avoid them we propose the following procedure. First, it is enough to get the linearized counterpart
\begin{equation}
A_0(U)\partial_t(\delta U)+\sum_{j=1}^3 A_j(U)\partial_j(\delta U)=0
\label{13}
\end{equation}
of system \eqref{first},  where $\delta U = (\delta p,\delta u,\delta H, \delta S)$ is the vector of perturbations and $U$ is now a {\it constant} vector, and $\rho =\rho (p,S)$, $e=e(\rho ,S)$, $v_i= u_i/\Gamma$, $\Gamma =(1-v^2)^{-1/2}$,  $1/a^2=\rho_p(\rho ,S)$ and $h=1+ e + (p/\rho )$ are constants. Below we write down system \eqref{13} for the fluid rest frame ($v=0$). After that we should only properly apply the Lorentz transformation and get system \eqref{13} in the LAB-frame.

The last equation in  \eqref{13} will be the linearization of \eqref{8},
\begin{equation}
\partial_t(\delta S) + (v,\nabla (\delta S))=0,
\label{14}
\end{equation}
and it is natural to believe that the matrices $A_{\alpha}$ in \eqref{13} will have the block structure
\begin{equation}
A_0=\begin{pmatrix} \widetilde{A}_0 & 0 \\ 0 & 1 \end{pmatrix},\quad A_j=\begin{pmatrix} \widetilde{A}_j & 0 \\ 0 & v_j \end{pmatrix}.
\label{14'}
\end{equation}
That is, our goal now is to find the symmetric matrices  $\widetilde{A}_{\alpha}$ of the system
\begin{equation}
\widetilde{A}_0(U)\partial_t(\delta V)+\sum_{j=1}^3 \widetilde{A}_j(U)\partial_j(\delta V)=0
\label{15}
\end{equation}
for $\delta V = (\delta p,\delta u,\delta H )$.

Let $U'=(p,0,H',S')= U|_{v=0}$ and $\delta V' = (\delta p',\delta u',\delta H') =\delta V|_{v=0}$. It is clear that the thermodynamical values coincide in the rest and LAB-frames: $p'=p$, $\rho' =\rho$, $h'=h$, $\delta p'=\delta p$, etc. We can easily write down the linearized system
\begin{align}
& \frac{1}{\rho a^2}\,\partial_t (\delta p')+{\rm div}\,(\delta u') =0,\label{16a}\\
& \bigl( \rho h+{H'}^2\bigr) \partial_t(\delta u')-H'(H',\partial_t(\delta u'))+\nabla (\delta p')+ H'\times (\nabla\times (\delta H'))=0,\label{16b}\\
& \partial_t(\delta H')- (H',\nabla)(\delta u')+H'{\rm div}\,(\delta u') =0
\label{16c}
\end{align}
associated with equations \eqref{2}, \eqref{3} and  \eqref{5} in the rest frame. Here we have taken into account \eqref{14} for $v=0$ and the divergence constraint $\div (\delta H')=0$. Equations \eqref{16a}--\eqref{16c} form the symmetric system
\begin{equation}
\widetilde{A}'_0(U')\partial_t(\delta V')+\sum_{j=1}^3\widetilde{A}'_j(U')\partial_j(\delta V')=0,\label{17}
\end{equation}
with
\begin{equation}
\widetilde{A}'_0=\left(
\begin{array}{ccc}
{\displaystyle\frac{1}{\rho a^2}} & 0 &0 \\[6pt]
0 & \mathcal{A}' &0 \\
0& 0& I\end{array}\right),\qquad \widetilde{A}'_j=\left(
\begin{array}{ccc}
0 & e_j^{\sf T} &0 \\[6pt]
e_j & 0 &\mathcal{K}_j^{\sf T} \\
0& \mathcal{K}_j& 0\end{array}\right),\label{17'}
\end{equation}
where
\[
\mathcal{A}'=\bigl(\rho h+{H'}^2\bigr)I-
H'\otimes H',\qquad
\mathcal{K}_j=-H'_jI+H'\otimes e_j
\]
and the magnetic field is transformed under the Lorentz transformation as
\begin{equation}
H'=\frac{H}{\Gamma} +\frac{\Gamma -1}{\Gamma v^2}(v,H)v.\label{18}
\end{equation}
The last formula is obtained by applying the Lorentz transformation to the 4-vector $(b^0,b)$:
\[
L\begin{pmatrix} (u,H) \\[3pt] \displaystyle{\frac{H}{\Gamma}}+(u,H)v \end{pmatrix} =
\begin{pmatrix} 0 \\[3pt] \displaystyle{\frac{H}{\Gamma} +\frac{\Gamma -1}{\Gamma v^2}(v,H)v} \end{pmatrix}= \begin{pmatrix} 0 \\ H' \end{pmatrix},
\]
where
\[
L=
\begin{pmatrix}
\Gamma & -\Gamma v\\[6pt]
-\Gamma v& I+{\displaystyle\frac{\Gamma -1}{v^2}\,v\otimes v}
\end{pmatrix}, \qquad
L^{-1}=
\begin{pmatrix}
\Gamma & \Gamma v\\[6pt]
\Gamma v& I+{\displaystyle\frac{\Gamma -1}{v^2}\,v\otimes v}
\end{pmatrix}.
\]

Now we should properly apply the Lorentz transformation to system \eqref{17}. The perturbations are transformed as
\begin{equation}
\delta V' =J\delta V  \label{19}
\end{equation}
and it is clear that the matrix $J$ has the form
\begin{equation}
J=\begin{pmatrix}
1 & 0 & 0\\
0&J_1& 0\\
0 & J_3& J_2
\end{pmatrix},\label{20}
\end{equation}
where $J_1$ is found by applying the Lorenz transformation to the perturbation $(\delta \Gamma,\delta u )=((v,\delta u), \delta u )$ of the 4-velocity $(\Gamma ,u)$, and $J_2$ and $J_3$ are found by the same procedure applied to the perturbation $(\delta b^0,\delta b)$ of the 4-vector $(b^0,b)$. Namely, we have:
\[
L\begin{pmatrix}
(v,\delta u) \\ \delta u
\end{pmatrix}
= \begin{pmatrix}
0 \\[6pt] \delta u +(v,\delta u){\displaystyle\frac{1-\Gamma }{\Gamma v^2}\,v}
\end{pmatrix} =
\begin{pmatrix}
0 \\ J_1\,\delta u
\end{pmatrix} =  \begin{pmatrix}
0 \\\delta u'
\end{pmatrix},
\]
\[
L\begin{pmatrix}
(u,\delta H) + (H,\delta u)\\ \delta b
\end{pmatrix}
=\begin{pmatrix}
{\displaystyle \frac{1}{\Gamma}\, (H,\delta u)}\\[6pt] J_3\,\delta u + J_2\,\delta H
\end{pmatrix} = \begin{pmatrix}
(H',\delta u')\\ \delta H'
\end{pmatrix},
\]
where
\begin{equation}
J_1= I + \frac{1-\Gamma }{\Gamma v^2}\,v\otimes v,\quad J_2= \frac{1}{\Gamma}\left( I+
\frac{\Gamma -1}{v^2}\,v\otimes v\right),\label{21}
\end{equation}
\[
J_3=-\frac{1}{\Gamma}\, H'\otimes v +(v,H')I + \frac{1-\Gamma }{\Gamma v^2}\,(v,H')v\otimes v ,
\]
and using \eqref{18}, we recalculate $J_3$:
\begin{equation}
J_3=-\frac{1}{\Gamma^2}\, H\otimes v +(v,H)I -(v,H)v\otimes v.\label{22}
\end{equation}

In view of \eqref{19},
\[
\widetilde{A}_{\alpha}=J^{\sf T}C_{\alpha}J,
\]
and for calculating the elements $c^{\alpha}_{kl}$ of the matrices $C_{\alpha}$ through the elements  $a^{\alpha}_{kl}$ of the matrices $\widetilde{A}'_{\alpha}$ we apply the Lorentz transformation to the 4-vector $(a^0_{kl},a^1_{kl},a^2_{kl},a^3_{kl})$:
\[
L^{-1}\begin{pmatrix} a^0_{kl} \\ \vdots \\ a^3_{kl}\end{pmatrix} =\begin{pmatrix} c^0_{kl} \\ \vdots \\ c^3_{kl}\end{pmatrix}.
\]
That is,
\begin{gather}
\widetilde{A}_{\alpha}=J^{\sf T}\Gamma \left(  \widetilde{A}'_0+\mathcal{G} \right) J =J^{\sf T}C_0J,\label{24}\\
\widetilde{A}_j=J^{\sf T}\left( \Gamma v_j\widetilde{A}'_0+\widetilde{A}'_j +\frac{\Gamma -1}{v^2}\,v_j \mathcal{G}\right) J =
v_j\widetilde{A}_0 +v_j\frac{1-\Gamma }{\Gamma v^2}\, J^{\sf T}\mathcal{G}J +J^{\sf T}\widetilde{A}'_jJ,\label{25}
\end{gather}
where
\[
\mathcal{G}= \sum_{k=1}^3v_k\widetilde{A}'_k=
\begin{pmatrix}
0& v^{\sf T} &0 \\
v & 0 &\mathcal{K}^{\sf T} \\
0& \mathcal{K} & 0
\end{pmatrix},\quad
C_0=\Gamma
\begin{pmatrix}
{\displaystyle\frac{1}{\rho a^2}} & v^{\sf T} &0 \\[6pt]
v & \mathcal{A}' &\mathcal{K}^{\sf T} \\
0& \mathcal{K} & I\end{pmatrix},\quad \mathcal{K}=H'\otimes v -(v,H')I.
\]

Using \eqref{18}, we recalculate:
\[
\mathcal{K}=\frac{1}{\Gamma}\, H\otimes v -(v,H)I +\frac{\Gamma -1}{\Gamma v^2}\,(v,H)v\otimes v,
\]
\[
\mathcal{K}_j =-\frac{1}{\Gamma}(H_jI -H\otimes e_j) + \frac{1-\Gamma }{\Gamma v^2}\,(v,H)
(v_j I-v\otimes e_j),
\]
\[
\mathcal{A}'= (\rho h +B^2)-\frac{1}{\Gamma^2}\,H\otimes H + \frac{\Gamma -1}{\Gamma^2v^2} \,(v,H)
\bigl( v\otimes H + H\otimes v\bigr) +\frac{(\Gamma -1)^2}{\Gamma^2 v^4}\,(v,H)^2v\otimes v.
\]
Then, from \eqref{14'}, \eqref{17'}, \eqref{20}--\eqref{25} after long calculations we find the symmetric matrices $A_{\alpha}$ in \eqref{27}.

As is known \cite{An}, natural physical restrictions guaranteeing the hyperbolicity of the RMHD system do not depend on the magnetic field and coincide with corresponding ones in relativistic gas dynamics. In our case, by direct calculations one can show that the hyperbolicity condition $A_0>0$ is equivalent to the condition $B_0>0$ for the relativistic Euler equations and holds if inequalities \eqref{hc} are satisfied (of course, by default we also assume that $|v|<1$). The last inequality in \eqref{9} is the {\it relativistic causality} condition.

\begin{remark}
Strictly speaking, to prove the equivalence of system \eqref{2}--\eqref{5} and \eqref{first} on smooth solutions we should also derive equations \eqref{2}--\eqref{5} from system \eqref{first}.
To this end we write down the subsystem for $H$ contained in \eqref{first} and deduce from this subsystem the divergence constraint \eqref{6} provided that it was satisfied for $t=0$.
The {remaining arguments are also very} similar to those for the non-relativistic case and we omit them.
\end{remark}

\begin{remark}
For interface problems for the RMHD equations when the interface moves with the velocity of plasma particles (e.g., for current-vortex sheets, see Sect. 4) it is useful to have the representation
\begin{equation}
A_j = v_jA_0 + G_j,\label{28}
\end{equation}
with
\[
G_j=
\begin{pmatrix}
0 & e_j^{\sf T}-v_jv^{\sf T} &0 \\[6pt]
e_j-v_jv & \mathcal{G}_j &{\mathcal{N}_j}^{\sf T} \\
0& \mathcal{N}_j & 0
\end{pmatrix},
\]
\begin{multline*}
\mathcal{G}_j= v_j\left\{ 2\frac{B^2}{\Gamma}\,v\otimes v -\frac{(v,H)}{\Gamma}
\bigl( v\otimes H + H\otimes v\bigr)\right\}\\
+\frac{H_j}{\Gamma}\left\{ \frac{1}{\Gamma^2}
\bigl( v\otimes H + H\otimes v\bigr) -2(v,H)(I-v\otimes v) \right\}
\\
+\frac{(v,H)}{\Gamma}\left( H\otimes e_j + e_j \otimes H\right) -\frac{B^2}{\Gamma}
\left( v\otimes e_j + e_j \otimes v\right).
\end{multline*}
\end{remark}

\section{Secondary symmetrization of the RMHD system}
\label{sec:3}

Now our goal is to find a different symmetrization of the RMHD equations. In some sense it will be a secondary symmetrization of system \eqref{first} which was already symmetric. This secondary symmetrization is a relativistic counterpart of the secondary symmetrization of the MHD system proposed in \cite{T05}. In Sect. 4 it will play a crucial role in finding a sufficient stability condition for relativistic current-vortex sheets.

We again consider the linear constant coefficient symmetric system \eqref{16a}--\eqref{16c}/\eqref{17} for the perturbation $\delta V'$ in the rest frame. For this system we have the standard conserved integral
\begin{equation}
\frac{\rm d}{{\rm d}t}\int_{\mathbb{R}^3}(\widetilde{A}'_0(U')\delta V',\delta V')\,{\rm d}x =0.\label{32}
\end{equation}
At the same time, we can derive another conserved integral. Indeed, multiplying \eqref{16b} and \eqref{16c} by $\delta H'$ and $\delta u'$ respectively, and using \eqref{16a} and the divergence constraint ${\rm div}\, \delta H' =0$, we get the following  conserved integral for $L_2$ solutions of the Cauchy problem:
\[
\frac{\rm d}{{\rm d}t}\int_{\mathbb{R}^3}\left\{ \frac{1}{c_s^2}(H',\delta u')\delta p -
(\mathcal{A}'\delta u', \delta H') \right\}{\rm d}x =0.
\]
The last integral can be rewritten as
\begin{equation}
\frac{\rm d}{{\rm d}t}\int_{\mathbb{R}^3}(\mathcal{L}_0\delta V',\delta V')\,{\rm d}x =0,
\label{33}
\end{equation}
with
\[
\mathcal{L}_0=
\begin{pmatrix}
0 & {\displaystyle \frac{{H'}^{\sf T}}{c_s^2}} &0 \\[6pt]
{\displaystyle \frac{H'}{c_s^2}} & 0 &-\mathcal{A}' \\[6pt]
0& -\mathcal{A}' & 0\end{pmatrix}.
\]

Combining \eqref{32} and \eqref{33}, we obtain the new conserved integral
\begin{equation}
\frac{\rm d}{{\rm d}t}\int_{\mathbb{R}^3}(\widetilde{\mathfrak{A}}'_0(U')\delta V',\delta V')\,{\rm d}x =0,
\label{34}
\end{equation}
where
\[
\widetilde{\mathfrak{A}}'_0=\widetilde{A}'_0+\lambda\mathcal{L}_0=
\begin{pmatrix}
{\displaystyle\frac{1}{\rho a^2}} & \lambda{\displaystyle \frac{{H'}^{\sf T}}{c_s^2}} &0 \\[6pt]
\lambda{\displaystyle \frac{H'}{c_s^2}} & \mathcal{A}' &-\lambda\mathcal{A}' \\[6pt]
0& -\lambda\mathcal{A}' & I\end{pmatrix},
\]
and $\lambda =\lambda (U')$ is an arbitrary constant.

Let us return to the nonlinear setting and consider the nonlinear counterpart
\begin{equation}
\widetilde{A}'_0(U')\partial_t V'+\sum_{j=1}^3\widetilde{A}'_j (U')\partial_j V'=0\label{35}
\end{equation}
of system \eqref{17}, where $V'=V|_{v=0}$ and $U'=U|_{v=0}=(V|_{v=0},S)$ are now the unknowns in the rest frame rather than constant vectors. For finding a symmetric system corresponding to the conserved integral \eqref{34}, where $\lambda =\lambda (U')$ is now a function of $U'$, we multiply \eqref{35} from the left by the matrix
\[
D'={D}'(U')=\widetilde{\mathfrak{A}}'_0(U'){(\widetilde{A}'_0(U'))}^{-1}=
\begin{pmatrix}
1 & {\displaystyle \frac{\lambda}{\rho a^2}\,{H'}^{\sf T}} &0 \\[9pt]
\lambda\rho h H' & I &-\lambda\mathcal{A}' \\[9pt]
0& -\lambda\mathcal{A}' & I\end{pmatrix}.
\]
{To get a symmetric system again, we 
add to the result of the} multiplication the vector $R'{\rm div}\, H'$, with $R'=R'(U')=
-\lambda (1,0,0,0, H')$.

Indeed,
\[
D'\widetilde{A}'_0=\widetilde{\mathfrak{A}}'_0={(\widetilde{\mathfrak{A}}'_0)}^{\sf T},\qquad
D'\widetilde{A}'_j=\begin{pmatrix}
{\displaystyle\frac{\lambda}{\rho a^2}\,H'_j} & e_j^{\sf T} &0 \\[9pt]
e_j & \lambda H'_j\mathcal{A}' &\mathcal{K}_j^{\sf T} \\[9pt]
-\lambda e_j& \mathcal{K}_j & -\lambda \mathcal{K}_j^{\sf T}
\end{pmatrix}.
\]
That is,
\begin{multline}
D'(U')\Bigl\{\widetilde{A}_0(U')\partial_tV'+\sum_{j=1}^3\widetilde{A}_j(U')\partial_jV'\Bigr\} +R'(U')\,{\rm div}\, H'\\
=\widetilde{\mathfrak{A}}'_0(U')\partial_tV'+\sum_{j=1}^3\widetilde{\mathfrak{A}}'_j(U')\partial_jV'=0,\label{36}
\end{multline}
with
\[
\widetilde{\mathfrak{A}}'_j(U')=\widetilde{A}'_j+\lambda\mathcal{L}_j=
\begin{pmatrix}
{\displaystyle\frac{\lambda}{\rho a^2}\,H'_j} & e_j^{\sf T} &-\lambda e_j^{\sf T} \\[6pt]
e_j & \lambda H'_j\mathcal{A}' &\mathcal{K}_j^{\sf T} \\
-\lambda e_j& \mathcal{K}_j& \lambda\mathcal{M}_j\end{pmatrix},
\]
and
\[
\mathcal{M}_j=H'_jI-(H'\otimes e_j+e_j\otimes H' ),\qquad \mathcal{L}_j=
\begin{pmatrix}
{\displaystyle\frac{H'_j}{\rho a^2}} & 0 &- e_j^{\sf T} \\[6pt]
0 &  H'_j\mathcal{A}' &0 \\
- e_j& 0& \mathcal{M}_j\end{pmatrix}.
\]

Note that
\[
D'>0\quad\Leftrightarrow\quad \lambda^2<\frac{c_s^2}{\rho a^2 + {H'}^2}\,,
\]
i.\ e., $\widetilde{\mathfrak{A}}'_0>0$ provided that (cf. \eqref{hc})
\begin{equation}
\rho>0,\quad p_{\rho} >0,\quad 0<c_s^2<1,\quad\mbox{and}\quad
\lambda^2<\frac{c_s^2}{\rho  a^2 + {H'}^2}\,.\label{37}
\end{equation}

To obtain the secondary symmetrization
\begin{equation}
\widetilde{\mathfrak{A}}_0(U)\partial_tV+\sum_{j=1}^3\widetilde{\mathfrak{A}}_j(U)\partial_jV=0\label{38}
\end{equation}
in the LAB-frame corresponding to \eqref{36} in the rest frame we again use the Lorenz transformation (see Sect. 2):
\begin{equation}
\widetilde{\mathfrak{A}}_0=J^{\sf T}\Gamma \Bigl( \widetilde{\mathfrak{A}}'_0+ \sum_{k=1}^3v_k\widetilde{\mathfrak{A}}'_k\Bigr) J=J^{\sf T} \bigl(C_0+\lambda\Gamma (\mathcal{L}_0+\mathcal{L})\,\bigr) J=
\widetilde{A}_0+\lambda \widetilde{B}_0=
J^{\sf T}\mathfrak{C}_0J,\label{39}
\end{equation}
\begin{multline}
\widetilde{\mathfrak{A}}_j=J^{\sf T}\Bigl( \Gamma v_j\widetilde{\mathfrak{A}}'_0+\widetilde{\mathfrak{A}}'_j +\frac{\Gamma -1}{v^2}\,v_j \sum_{k=1}^3v_k\widetilde{\mathfrak{A}}'_k\Bigr) J\\
=\widetilde{A}_j + v_j\lambda \widetilde{B}_0 +
\lambda\Bigl( v_j\frac{1-\Gamma }{\Gamma v^2}\, J^{\sf T}\mathcal{L}J +J^{\sf T}\mathcal{L}_jJ\Bigr),\label{40}
\end{multline}
where $\lambda =\lambda (U)$,
\[
\mathcal{L}=\sum_{k=1}^3v_k\mathcal{L}_k=
\begin{pmatrix}
{\displaystyle\frac{(v,H')}{\rho a^2}} & 0 &- v^{\sf T} \\[6pt]
0 &  (v,H')\mathcal{A}' &0 \\
- v& 0& \mathcal{N}\end{pmatrix},\quad \mathcal{N}=(v,H')I-(H'\otimes v+v\otimes H' ),
\]
\[
\mathfrak{C}_0=C_0+\lambda\Gamma (\mathcal{L}_0+\mathcal{L})=
\Gamma
\begin{pmatrix}
{\displaystyle\frac{1+\lambda (v,H')}{\rho a^2}} & {\displaystyle v^{\sf T}+\lambda\frac{{H'}^{\sf T}}{c_s^2}} &- \lambda v^{\sf T} \\[6pt]
{\displaystyle v+\lambda\frac{{H'}}{c_s^2}} &  \bigl(1+\lambda (v,H')\bigr)\mathcal{A}' &\mathcal{K}^{\sf T}-\lambda\mathcal{A}' \\
-\lambda v& \mathcal{K}-\lambda\mathcal{A}'& I+\lambda\mathcal{N}\end{pmatrix}.
\]

Note that $(v,H') =(v,H)$ and using \eqref{18} all the other values appearing in the formulas above can be rewritten in terms of $H$. System \eqref{38} together with equation \eqref{8} form the symmetric system
\begin{equation}
\mathfrak{A}_0(U)\partial_tU+\sum_{j=1}^3\mathfrak{A}_j(U)\partial_jU=0\label{41}
\end{equation}
with
\begin{equation}
\mathfrak{A}_0=\begin{pmatrix} \widetilde{\mathfrak{A}}_0 & 0 \\ 0 & 1 \end{pmatrix},\quad \mathfrak{A}_j=\begin{pmatrix} \widetilde{\mathfrak{A}}_j & 0 \\ 0 & v_j \end{pmatrix},
\label{42}
\end{equation}
where, if necessary, the symmetric matrices $\widetilde{\mathfrak{A}}_{\alpha}$ can be written in a concrete form (in terms of dyadic products) like matrices $A_{\alpha}$ in  \eqref{first}. System \eqref{41} is nothing else than our secondary symmetrization \eqref{second} in Theorem \ref{t2}.

Regarding the hyperbolicity condition $\mathfrak{A}_0>0$ for system \eqref{41}, it is equivalent to the requirement $\mathfrak{C}_0>0$ which holds provided that
\begin{equation}
\rho>0,\quad p_{\rho} >0,\quad 0<c_s^2<1,\quad\mbox{and}\quad
\lambda^2<\frac{c_s^2}{\rho  a^2 + B^2}\label{43}
\end{equation}
(we drop technical calculations). It is natural that the hyperbolicity condition \eqref{43} for the secondary symmetrization \eqref{41} in the LAB-frame coincides with the hyperbolicity condition \eqref{37} in the rest frame because ${H'}^2=B^2$.

\section{Stability of relativistic current-vortex sheets}
\label{sec:4}

The initial boundary value problem (in fact, the free boundary problem) for relativistic current-vortex sheets is formulated as follows. We consider the RMHD equations for $t\in [0,T]$ in the unbounded space domain $\mathbb{R}^3$ and assume that
\begin{equation}
\Sigma (t)=\{ x^1=\varphi (t,x')\}\label{44}
\end{equation}
is a smooth hypersurface in $[0,T]\times\mathbb{R}^3$, where
$x'=(x^2,x^3)$ are tangential coordinates. We assume that $\Sigma (t)$ is a surface of strong discontinuity for the conservation laws \eqref{2}--\eqref{5},
i.\ e., we are interested in solutions of \eqref{2}--\eqref{5} that are smooth on either side of $\Sigma (t)$. To be weak solutions of \eqref{2}--\eqref{5} such piecewise smooth solutions should satisfy Rankine-Hugoniot jump conditions. Since we are interested only in current-vortex sheets, we do not write down here the
RMHD Rankine-Hugoniot conditions in a general case which covers, in particular, shock waves.

As in non-relativistic settings, for relativistic current-vortex sheets we require that
the hypersurface $\Sigma (t)$ moves with the velocity of plasma particles at the boundary and the magnetic field on $\Sigma (t)$ is parallel to $\Sigma (t)$. With these requirements the Rankine-Hugoniot conditions give the boundary conditions
\begin{align}
& \partial_t\varphi =v^{\pm}_N, \label{45}\\
& [q] =0, \label{46}\\
& H_{\rm N}^{\pm}=0\qquad \mbox{on}\ \Sigma (t),\label{47}
\end{align}
where $[g]=g^+|_{\Sigma}-g^-|_{\Sigma}$ denotes the jump of $g$, with $g^{\pm}:=g$ in $\Omega^{\pm}(t)=\{ x^1\gtrless \varphi(t,{x}')\}$; $v^{\pm}_N=(v^{\pm},N)$, $H^{\pm}_N=(H^{\pm},N)$, and $N=(1,-\partial_2\varphi ,-\partial_3\varphi )$. Again, as
{in non-relativistic settings} \cite{T09}, one can show that \eqref{47} are not real boundary conditions and should be considered as restrictions on the initial data
\begin{equation}
{U}^{\pm} (0,{x})={U}_0^{\pm}({x}),\quad {x}\in \Omega^{\pm} (0),\qquad \varphi (0,{x}')=\varphi _0({x}'),\quad {x}'\in\mathbb{R}^2.\label{48}
\end{equation}

\begin{remark}
Our assumption that the hypersurface \eqref{44} has the form of a graph is not a strong restriction for special relativity. Regarding the case of general relativity, $\Sigma$ should be a compact  codimension-1 surface. However, as for shock waves we can follow Majda's arguments \cite{Maj83b} for extending the results to a compact hypersurface. On the other hand, even for general relativity, this assumption is
{still not} a strong restriction because, as in \cite{FR}, we can consider the interface not only locally in time (in the sense of local-in-time existence of piecewise smooth solutions), but also locally in space (for compactly supported initial data).
\end{remark}

As in the non-relativistic case \cite{T05,T09},  the linear stability condition for planar current-vortex sheets which should be satisfied at each point of the initial hypersurface $\Sigma(0)$ is the basic assumption providing the local-in-time existence of smooth (nonplanar) current-vortex sheet solutions,
i.\ e., the existence of a solution $(U^+,U^- ,\varphi )$ of the free boundary value problem \eqref{2}--\eqref{5}, \eqref{45}--\eqref{48},  where $U^{\pm}:=U$ in $\Omega^{\pm}(t)$, and ${U}^{\pm}$ is smooth in $\Omega^{\pm}(t)$. As for non-relativistic current-vortex sheets \cite{T05}, we formulate the constant coefficient linearized stability problem associated to the original nonlinear problem \eqref{2}--\eqref{5}, \eqref{45}--\eqref{48}:
\begin{equation}
\widehat{A}^{\pm}_0\partial_t{U}^{\pm}+\sum_{j=1}^{3}\widehat{A}^{\pm}_j\partial_j{U}^{\pm}=0 \qquad \mbox{if}\ x\in\mathbb{R}^3_{\pm},\label{49}
\end{equation}
\begin{equation}
\partial_t\varphi={v}_1^{\pm}-\hat{v}_2\partial_2\varphi-\hat{v}_3^{\pm}\partial_3\varphi , \quad [q]=0,   \qquad \mbox{if}\ x^1=0,
\label{50}
\end{equation}
\begin{equation}
{U}^{\pm} (0,{x})={U}_0^{\pm}({x}),\quad {x}\in \mathbb{R}^3_{\pm},\qquad \varphi (0,{x}')=\varphi _0({x}'),\quad {x}'\in\mathbb{R}^2,\label{51}
\end{equation}
where $(U^+,U^-,\varphi )$ is now a perturbation of the constant current-vortex sheet solution $(\widehat{U}^+,\widehat{U}^-,0)$ corresponding to a planar relativistic current-vortex sheet with the equation $x^1=0$ (all the hat values are associated with this solution and without loss of generality we take $\hat{\varphi}=0$), $\mathbb{R}^3_{\pm}=\{\pm x^1> 0,\ x'\in\mathbb{R}^2\}$, $\widehat{A}_{\alpha}={A}_{\alpha}(\widehat{U})$, $[q]=q^+_{|x^1=0}-q^-_{|x^1=0}$. Since $(\widehat{U}^+,\widehat{U}^-,0)$ is assumed to be a current-vortex sheet solution, in view of \eqref{45} and \eqref{47},
\[
U^{\pm}=(\hat{p}^{\pm},0,\hat{v}_2^{\pm},\hat{v}_3^{\pm},0,\widehat{H}_2^{\pm},\widehat{H}_3^{\pm},
\widehat{S}).
\]
Again, as in non-relativistic settings \cite{T05,T09}, we can show that the solution of problem \eqref{49}-\eqref{51}
{satisfies} the linearized conditions \eqref{47},
\begin{equation}
{H}_1^{\pm}=\widehat{H}_2\partial_2\varphi+\widehat{H}_3^{\pm}\partial_3\varphi    \qquad \mbox{if}\ x^1=0,
\label{52}
\end{equation}
provided that these conditions were satisfied by the initial data \eqref{51}.

The linearized stability of a planar relativistic current-vortex sheet means the well-posedness of problem \eqref{49}-\eqref{51} and is
{equivalent to} the fulfilment of the Kreiss-Lopatinski condition (see, e.g., \cite{FR,Maj83b}). As was shown in \cite{T05}, for non-relativistic current-vortex sheets the Kreiss-Lopatinski condition can be satisfied only in a weak sense,
i.\ e., they can be only {\it neutrally} stable. Here we do not formally show that the same is true in the relativistic case, but it is natural to expect that this is really so,
i.\ e., the uniform Kreiss-Lopatinski condition is never satisfied for relativistic current-vortex sheets. Due to technical reasons \cite{T05} it is almost impossible to study the linear stability of current-vortex sheets by the standard method of normal modes. It is clear that for the relativistic case the situation is even worse. Therefore, as in \cite{T05}, we will use the energy method for finding a sufficient neutral stability condition for relativistic current-vortex sheets.

In the application of the energy method the secondary symmetrization \eqref{41} found in Sect. 3 plays the crucial role. The linearization of \eqref{41}
{about} the constant solution reads:
\begin{equation}
\widehat{\mathfrak{A}}^{\pm}_0\partial_tU^{\pm}+\sum_{j=1}^3\widehat{\mathfrak{A}}^{\pm}_j\partial_jU^{\pm}=0
\qquad \mbox{if}\ x\in\mathbb{R}^3_{\pm},
\label{53}
\end{equation}
where $\widehat{\mathfrak{A}}^{\pm}_{\alpha}=\mathfrak{A}_{\alpha}(\widehat{U}^{\pm})$. For systems \eqref{53} standard arguments of the energy method give the {\it energy integral identity}
\begin{equation}
\frac{{\rm d}}{{\rm d} t} I(t)-\int_{{\mathbb R}^2}
\bigl[(\widehat{\mathfrak{B}}_1U ,U )\bigr]\bigl|_{x^1=0}\,{\rm d} x'=0. \label{54}
\end{equation}
where
\[
\bigl[(\widehat{\mathfrak{B}}_1U ,U )\bigr]\bigl|_{x^1=0} =
(\widehat{\mathfrak{B}}_1^+U^+ ,U^+ )|_{x^1=0}-(\widehat{\mathfrak{B}}_1^-U^- ,U ^- )|_{x^1=0},\quad
\widehat{\mathfrak{B}}_1^{\pm}=\beta^{\pm}\widehat{\mathfrak{A}}_1^{\pm}
\]
\[
I(t)=\int_{{\mathbb R}^3_+}\beta^+(\widehat{\mathfrak{A}}_0^+ U^+ ,U^+ )\,{\rm d}{x}
+ \int_{{\mathbb R}^3_-}\beta^-(\widehat{\mathfrak{A}}_{0}^- U^- ,U^- )\,{\rm d}{x},
\]
and the constants $\beta^+>0$ and  $\beta^->0$ will be chosen later on.
As follows from \eqref{54}, the {\it dissipativity condition}
$
\bigl[(\widehat{\mathfrak{B}}_1U ,U )\bigr]\bigl|_{x^1=0}\leq 0
$
satisfied for an unperturbed flow $(\widehat{U}^+,\widehat{U}^-)$ guarantees the stability of the corresponding planar current-vortex sheet. Following \cite{T05,T09}, we can choose such $\lambda^{\pm}=\lambda (U^{\pm})$ in the secondary symmetrization \eqref{41} that
\begin{equation}
\bigl[(\widehat{\mathfrak{B}}_1U ,U )\bigr]\bigl|_{x^1=0}=0,\label{55}
\end{equation}
and the last necessary inequality in the hyperbolicity conditions \eqref{43} for chosen $\lambda^{\pm}$ gives us a sufficient stability condition.

Before choosing $\lambda^{\pm}$ we have to compute the boundary matrices $\widehat{\mathfrak{A}}_1^{\pm}$. Actually, we {do not need
these matrices themselves, but} the quadratic forms $(\widehat{\mathfrak{A}}_1^{\pm}U^{\pm} ,U^{\pm} )|_{x^1=0}$. Let us compute them even for a general basic state  $(\widehat{U}^+(t,x),\widehat{U}^-(t,x),\hat{\varphi}(t,x'))$ satisfying conditions \eqref{45} and \eqref{47}. Then the boundary matrices read
\[
\widehat{\mathfrak{A}}^{\pm}_N= \widehat{\mathfrak{A}}^{\pm}_1-\widehat{\mathfrak{A}}^{\pm}_0\partial_t\hat{\varphi}-\widehat{\mathfrak{A}}^{\pm}_2\partial_2\hat{\varphi}-
\widehat{\mathfrak{A}}^{\pm}_3\partial_3\hat{\varphi}.
\]
We first calculate the quadratic forms $(\widehat{\mathfrak{A}}_N^{\pm}U^{\pm} ,U^{\pm} )|_{x^1=0}$ for $\lambda^{\pm}=0$,
i.\ e., the quadratic forms $(\widehat{A}_N^{\pm}U^{\pm} ,U^{\pm} )|_{x^1=0}$.

Taking into account \eqref{45}, \eqref{47} and representation \eqref{28}, we have
\[
\widehat{A}_N^{\pm}|_{x^1=0}= \widehat{G}_N^{\pm}|_{x^1=0}=\Bigl({G}_1(\widehat{U}^{\pm})-{G}_2(\widehat{U}^{\pm})\partial_2\hat{\varphi}
-{G}_3(\widehat{U}^{\pm})\partial_3\hat{\varphi}\Bigr)\Bigr|_{x^1=0}.
\]
After some algebra we find
\[
\bigl(\widehat{A}_N^{\pm}U^{\pm},U^{\pm}\bigr)\bigr|_{x^1=0}=2\widehat{\Gamma}^{\pm}q^{\pm}v_N^{\pm}|_{x^1=0},
\]
where $v_N^{\pm}=(v^{\pm},\widehat{N})$, $\widehat{N}=(1,-\partial_2\hat{\varphi},-\partial_3\hat{\varphi})$,
\[
 q^{\pm}= p^{\pm} +\frac{1}{2}\delta (B^{\pm})^2=p^{\pm}+\frac{1}{\widehat{\Gamma}^{\pm}}
\left( (\hat{b}^{\pm},H^{\pm}) +(\hat{v}^{\pm},\widehat{H}^{\pm})(\widehat{H}^{\pm},u^{\pm})-(\widehat{B}^{\pm})^2(\hat{v}^{\pm},u)\right),
\]
and $\delta (B^{\pm})^2$ is the perturbation of $B^2$ for $\pm x^1>0$.

Now we consider the case $\lambda^{\pm}\neq 0$:
\[
\bigl(\widehat{\mathfrak{A}}_N^{\pm}U,U\bigr)\bigr|_{x^1=0}=
\bigl(\widehat{{A}}_N^{\pm}U,U\bigr)\bigr|_{x^1=0}+\lambda^{\pm}
\bigl(\widehat{\mathfrak{M}}^{\pm}U,U\bigr)\bigr|_{x^1=0},
\]
where
\[
\bigl(\widehat{\mathfrak{M}}^{\pm}U,U\bigr)\bigr|_{x^1=0}=
\left.
\left(\Bigl\{\hat{\varphi}_t\frac{1-\widehat{\Gamma}^{\pm} }{\widehat{\Gamma}^{\pm}(\hat{v}^{\pm})^2}\,\mathcal{L}(\widehat{U}^{\pm})+ \widehat{\mathcal{L}}_N^{\pm} \Bigr\}J(\widehat{U}^{\pm})\,U,J(\widehat{U}^{\pm})\,U\right)\right|_{x^1=0},
\]
\[
\widehat{\mathcal{L}}_N^{\pm}=\mathcal{L}_1(\widehat{U}^{\pm})-\mathcal{L}_2(\widehat{U}^{\pm})\partial_2\hat{\varphi}-
\mathcal{L}_3(\widehat{U}^{\pm})\partial_3\hat{\varphi}.
\]
Omitting (not extremely long) calculations, we obtain
\begin{multline*}
\bigl(\widehat{\mathfrak{M}}^{\pm}U,U\bigr)\bigr|_{x^1=0}=
-2\bigl(p^{\pm} +((\widehat{H}')^{\pm},(H')^{\pm})\bigr)\left.\left(((H')^{\pm},\widehat{N})+\hat{\varphi}_t\frac{1-\widehat{\Gamma}^{\pm} }{\widehat{\Gamma}^{\pm}(\hat{v}^{\pm})^2}\,(\hat{v}^{\pm},(H')^{\pm})\right)\right|_{x^1=0}
\\[6pt]
=-2\widehat{\Gamma}^{\pm}q^{\pm}\left.\left(\frac{H_N^{\pm}}{{\bigl(\widehat{\Gamma}^{\pm}\bigr)}^2}+(\hat{v}^{\pm},\widehat{H}^{\pm})v_N^{\pm}
\right)\right|_{x^1=0}.
\end{multline*}
That is,
\[
\bigl(\widehat{\mathfrak{A}}_N^{\pm}U,U\bigr)\bigr|_{x^1=0}=
2\widehat{\Gamma}^{\pm}q^{\pm}\left.\left(\bigl(1-\lambda^{\pm} (\hat{v}^{\pm},\widehat{H}^{\pm})\bigr)\,v_N^{\pm}-\frac{\lambda^{\pm}}{(\widehat{\Gamma}^{\pm})^2}\,H_N^{\pm}\right)\right|_{x^1=0},
\]
and for the planar discontinuity with $\hat{\varphi}=0$
\begin{equation}
\bigl(\widehat{\mathfrak{A}}_1^{\pm}U,U\bigr)\bigr|_{x^1=0}=
2\widehat{\Gamma}^{\pm}q^{\pm}\left.\left(\bigl(1-\lambda^{\pm} (\hat{v}^{\pm},\widehat{H}^{\pm})\bigr)\,v_1^{\pm}-\frac{\lambda^{\pm}}{(\widehat{\Gamma}^{\pm})^2}\,H_1^{\pm}\right)\right|_{x^1=0}.
\label{56}
\end{equation}

We now choose $\beta^{\pm}$ in \eqref{54}:
\[
\beta^{\pm}=\frac{1}{\widehat{\Gamma}^{\pm}\bigl(1-\lambda^{\pm} (\hat{v}^{\pm},\widehat{H}^{\pm})\bigr)}.
\]
We assume that the hyperbolicity condition \eqref{43} holds for $U=\widehat{U}^{\pm}$, i.\ e., $\widehat{\mathfrak{A}}_0^{\pm}>0$. The last inequality in \eqref{43} guarantees that $1-\lambda^{\pm} (\hat{v}^{\pm},\widehat{H}^{\pm})>0$,
i.\ e., $\beta^{\pm}>0$. Then $I(t)>0$ and it follows from \eqref{56} that
\[
\bigl[(\widehat{\mathfrak{B}}_1U ,U )\bigr]\bigl|_{x^1=0}=2q^+[v_1-\tilde{\lambda}H_1]|_{x^1=0}.
\]
where
\begin{equation}
\tilde{\lambda}^{\pm}=\frac{\lambda^{\pm}}{(\widehat{\Gamma}^{\pm})^2\bigl(1-\lambda^{\pm} (\hat{v}^{\pm},\widehat{H}^{\pm})\bigr)}.\label{lamb}
\end{equation}

The boundary conditions \eqref{50} and \eqref{52} imply
\[
[v_1-\tilde{\lambda}H_1]|_{x^1=0}= \bigl(\bigl[\hat{v}^{\|}-\tilde{\lambda}\widehat{H}^{\|}\bigr],\nabla_{\rm tan}\varphi \bigr),
\]
where $\hat{v}^{\|\pm}=(\hat{v}_2^{\pm},\hat{v}_3^{\pm})$, $\widehat{H}^{\|\pm}=(\widehat{H}_2^{\pm},\widehat{H}_3^{\pm})$ and $\nabla_{\rm tan}=(\partial_2,\partial_3)$. As for non-relativistic current-vortex sheets \cite{T05,T09}, we choose $\tilde{\lambda}^{\pm}$ such that
\begin{equation}
\bigl[\hat{v}^{\|}-\tilde{\lambda}\widehat{H}^{\|}\bigr]=0,\label{57}
\end{equation}
i.\ e., \eqref{55} holds. This gives us an a priori $L_2$ estimate if the chosen $\tilde{\lambda}^{\pm}$ satisfy the last inequality in \eqref{43}. At the same time, this inequality for the chosen $\tilde{\lambda}^{\pm}$ gives a {\it sufficient} stability condition for planar relativistic current-vortex sheets.

Assume that the tangential magnetic fields $\widehat{H}^{\|\pm}$ are nonzero and nonparallel to each other:
\begin{equation}
\widehat{H}_2^+\widehat{H}_3^- - \widehat{H}_3^+\widehat{H}_2^- \neq 0.\label{58}
\end{equation}
For non-relativistic current-vortex sheets the violation of \eqref{58} corresponds to the transition to instability. It is natural to expect that the same is true in the relativistic case. It follows from \eqref{57} and \eqref{58} that
\[
|\tilde{\lambda}^{\pm}|=\frac{|[\hat{v}^{\|}]|\,|\sin\varphi^{\mp}|}{|\widehat{H}^{\|\pm}|\,|\sin (\varphi^+ -\varphi^-)|}\,,
\]
where
\begin{equation}
\label{not}
{v}^{\|\pm}=(v_2^{\pm},v_3^{\pm}),\quad
\widehat{H}^{\|\pm}=(\widehat{H}_2^{\pm},\widehat{H}_3^{\pm}),\quad [\hat{v}^{\|}]=\hat{v}^{\|+}-\hat{v}^{\|-},\quad
\cos\varphi^{\pm}=\frac{([\hat{v}^{\|}],\widehat{H}^{\|\pm})}{ |[\hat{v}^{\|}]|\,|\widehat{H}^{\|\pm}|}\,.
\end{equation}

In view of \eqref{lamb},
\[
\lambda^{\pm}=\frac{\hat{\lambda}^{\pm}}{1+\hat{\lambda}^{\pm} (\hat{v}^{\pm},\widehat{H}^{\pm})},\qquad\mbox{with}\quad \hat{\lambda}^{\pm}=(\widehat{\Gamma}^{\pm})^2\tilde{\lambda}^{\pm}.
\]
The last condition in \eqref{43} for $\lambda^{\pm}$ gives
\[
|\lambda^{\pm}|< m^{\pm},\qquad\mbox{with}\quad m^{\pm}=\frac{\hat{c}^{\pm}_s}{\sqrt{\hat{\rho}^{\pm}  (\hat{a}^{\pm})^2 + (\widehat{B}^{\pm})^2}}\,.
\]
We can show that in terms of $\hat{\lambda}^{\pm}$ the last inequalities read
\[
|\hat{\lambda}^{\pm}|<\frac{1}{{\displaystyle\frac{1}{m^{\pm}}}+|(\hat{v}^{\pm},\widehat{H}^{\pm})|}\,,
\]
i.\ e.,
\[
|\tilde{\lambda}^{\pm}|<\frac{1/(\widehat{\Gamma}^{\pm})^2}{{\displaystyle\frac{1}{m^{\pm}}}+|(\hat{v}^{\pm},\widehat{H}^{\pm})|}\,.
\]
Then, we can finally write down the sufficient stability condition in the same form as for the non-relativistic current-vortex sheets (we drop the hats,
i.\ e., write down the condition to be satisfied for the initial data in a counterpart of the nonlinear existence theorem from \cite{T09}):
\begin{equation}
G(U^+,U^-)>0, \label{sc}
\end{equation}
where
\[
G(U^+,U^-)=|\sin(\varphi^+-\varphi^-)|
\min\left\{\frac{\gamma^+}{|\sin{\varphi^-}|}\,,\,
\frac{\gamma^-}{|\sin{\varphi^+}|}\right\}-|[{v}^{\|}]|,
\]
but now (in the {\it relativistic} case)
\[
\gamma^{\pm}=\frac{|{H}^{\|\pm}|\,c_s^{\pm}\bigl(1-(v^{\pm})^2\bigr)}{\sqrt{\rho^{\pm}  (a^{\pm})^2 + ({B}^{\pm})^2}+c_s^{\pm}\,|(\hat{v}^{\pm},\widehat{H}^{\pm})|}.
\]
Recall that we also assume (cf. \eqref{58})
\begin{equation}
H_2^+H_3^- - H_3^+H_2^- \geq \epsilon >0,
\label{sc'}
\end{equation}
where $\epsilon$ is a fixed constant.

If in $\gamma^{\pm}$ we formally set $v^{\pm}=0$  and the specific enthalpies $h^{\pm}=1$ (this corresponds to the non-relativistic limit), then \eqref{sc} coincides with the sufficient stability condition for the classical current-vortex sheets \cite{T05,T09}.

Since current-vortex sheets can be only neutrally stable, there appears a loss of derivatives phenomenon in a priori estimates. Therefore, the nonlinear existence theorem was proved in \cite{T09} by Nash-Moser iterations. In the relativistic case, all the arguments towards the proof of the
local-in-time existence of smooth current-vortex sheet solutions, provided that the stability condition \eqref{sc}, \eqref{sc'} is satisfied at each point of the initial discontinuity, are absolutely the same as in the non-relativistic case in \cite{T09}.

\end{document}